\documentclass[12pt]{amsart}
\usepackage{amsmath,amsfonts,amssymb,mathrsfs}
\usepackage{times,array,delarray,graphicx,dsfont,xspace,a4wide}
\usepackage[numbers,sort&compress]{natbib}
\usepackage[colorlinks=true,linkcolor=blue,urlcolor=red,citecolor=blue,
  pdfauthor=Matthew\ Macauley\ and\ Jon\ McCammond\ and\ Henning\ S.\
  Mortveit, pdftitle=Order Independence in Asynchronous
  Cellular Automata, pdfsubject=Discrete\
  Mathematics,pdfkeywords=SDS]{hyperref}

\theoremstyle{definition}
\newtheorem{defn}{Definition}[section]
\newtheorem{prop}[defn]{Proposition}

\newtheorem{thm}[defn]{Theorem}

\newtheorem{rem}[defn]{Remark}
\newtheorem{exmp}[defn]{Example}
\numberwithin{equation}{section}

\newcommand{\y}{\mathbf{y}}

\newcommand{\zer}{\mathbf{0}}
\newcommand{\one}{\mathbf{1}}
\newcommand{\F}{\mathbb{F}}

\newcommand{\real}{\mathbb{R}}
\renewcommand{\frak}{\mathfrak}
\newcommand{\boxsp}{\hspace{0.26cm}}

\newcommand{\inv}{\mathsf{inv}}

\newcommand{\I}{I}
\newcommand{\refl}{\mathsf{refl}}

\newcommand{\R}{R}

\newcommand{\per}{{\sf Per}\xspace }
\newcommand{\ca}{{\sf CA}\xspace}
\newcommand{\cas}{{\sf CA}s\xspace}
\newcommand{\sds}{{\sf SDS}\xspace}
\newcommand{\sdss}{{\sf SDS}s\xspace}
\newcommand{\aca}{{\sf ACA}\xspace}
\newcommand{\acas}{{\sf ACA}s\xspace}
\newcommand{\Circ}{{\sf Circ}\xspace }
\newcommand{\wolf}{{\sf wolf}\xspace }
\newcommand{\Wolf}{{\sf Wolf}\xspace }
\newcommand{\WOLF}{{\sf \frak{Wolf}}\xspace }
\newcommand{\Fy}{\mathfrak{F}_Y}

\begin{document}
\bibliographystyle{amsplain}

\title[Order Independence in \acas]{Order Independence in Asynchronous
  Cellular Automata}

\author{
M.~Macauley\qquad\qquad J.~McCammond
\qquad\qquad H.S.~Mortveit
}

\date{\today}

\begin{abstract}
A \emph{sequential dynamical system}, or \sds, consists of an
undirected graph $Y$, a vertex-indexed list of local functions $\Fy$,
and a permutation $\pi$ of the vertex set (or more generally, a word
$\omega$ over the vertex set) that describes the order in which these
local functions are to be applied.  In this article we investigate the
special case where $Y$ is a circular graph with $n$ vertices and all
of the local functions are identical.  The $256$ possible local
functions are known as \emph{Wolfram rules} and the resulting
sequential dynamical systems are called \emph{finite asynchronous
  elementary cellular automata}, or \acas, since they resemble
classical elementary cellular automata, but with the important
distinction that the vertex functions are applied sequentially rather
than in parallel.  An \aca is said to be \emph{$\pi$-independent} if
the set of periodic states does not depend on the choice of $\pi$, and
our main result is that for all $n>3$ exactly $104$ of the $256$
Wolfram rules give rise to a $\pi$-independent \aca.  In 2005 Hansson,
Mortveit and Reidys classified the $11$ symmetric Wolfram rules with
this property.  In addition to reproving and extending this earlier
result, our proofs of $\pi$-independence also provide significant
insight into the dynamics of these systems.
\end{abstract}

\keywords{Cellular automata, periodic points, potential function,
  sequential dynamical systems, update order, Wolfram rules.}

\subjclass[2000]{37B99,68Q80}

\maketitle

Our main result, as recorded in Theorem~\ref{thm:104}, is a complete
classification of the Wolfram rules that for all $n>3$ lead to a
$\pi$-independent finite asynchronous elementary cellular automaton,
or \aca.  The structure of the article is relatively straightforward.
The first two sections briefly describe how an \aca can be viewed as
either a special type of sequential dynamical system or as a modified
version of a classical elementary cellular automaton. These two
sections also contain the background definitions and notations needed
to carefully state our main result.  Next, we introduce several new
notations for Wolfram rules in order to make certain patterns easier
to discern, and we significantly reduce the number of cases we need to
consider by invoking the notion of dynamical equivalence.
Sections~\ref{sec:major} and~\ref{sec:minor} contain the heart of the
proof. The former covers four large classes of rules whose members are
$\pi$-independent for similar reasons, and the latter finishes off
three pairs of unusual cases that exhibit more intricate behavior
requiring more delicate proofs. The final section contains remarks
about directions for future research.

\section{Sequential Dynamical Systems}\label{sec:sds}

Cellular automata, or \cas, are discrete dynamical systems that
have been thoroughly studied by both professional and amateur
mathematicians.%
\footnote{Stanis{\l}aw Ulam and John von Neumann were the first to study
  such systems, which they did while working at Los Alamos National
  Laboratory in the 1940s \cite{Neumann:66a}. The German computer
  scientist Konrad Zuse proposed in 1969 that the universe is
  essentially one big cellular automaton~\cite{Zuse:69}. In the 1970s,
  John Conway invented the Game of Life, a two-dimensional \ca, that
  was later popularized by Martin Gardner~\cite{Gardner:70}. Beginning
  in 1983, Stephen Wolfram published a series of papers devoted to
  developing a theory of \cas and their role in science
  \cite{Martin:84,Wolfram:83,Wolfram:84,Wolfram:86a}. This is also a
  central theme in Wolfram's $1280$-page book \emph{A New Kind of
    Science}, published in 2002.}
They are defined over regular grids of cells and each cell can take on
one of a finite number of states. In addition, each cell has an
\emph{update rule} that takes its own state and the states of its
neighbors as input, and at each discrete time step, the rules are
applied and all vertex states are simultaneously updated.

In the late 1990s, a group of scientists at Los Alamos invented a new
type of discrete dynamical system that they called \emph{sequential
  dynamical systems}, or \sdss.  In an \sds the regular grids used to
define cellular automata are replaced by arbitrary undirected graphs,
and the local functions are applied sequentially rather than in
parallel.  Their initial motivation was to develop a mathematical
foundation for the analysis, simulation and implementation of various
socio-technological systems~\cite{Barrett:99a, Barrett:00a,
  Mortveit:01b, Barrett:03a}.

An \sds has three components: an undirected graph $Y$, a list of local
functions $\Fy$, and an update order $\omega$.  Start with a
simple undirected graph $Y$ with $n$ vertices, label the vertices from
$1$ to $n$, and recall that the \emph{neighbors} of a vertex are those
vertices connected to it by an edge.
If $\F$ is a finite field and every vertex is assigned a value from
$\F$, then a global state of the system is described by an $n$-tuple
$\y$ whose $i^{\rm th}$ coordinate indicates the current state of the
vertex $i$.  The set of all possible states is the vector space
$\F^n$.

\begin{defn}[Local functions]
A function $F\colon\F^n\to \F^n$ is called \emph{$Y$-local at $i$} if
(1) for each $\y\in \F^n$, $F(\y)$ only alters the $i^{\rm th}$
coordinate of $\y$ and (2) the new value of the $i^{\rm th}$
coordinate only depends on the coordinates of $\y$ corresponding to
$i$ and its neighbors in $Y$.  Other names for such a function include
\emph{local function} and \emph{update rule}.  We use $\Fy$ to
denote a list of local functions that includes one for each vertex of
$Y$.  More precisely, $\Fy = (F_1, F_2, \ldots, F_n)$ where
$F_i$ is a function that is $Y$-local at $i$.
\end{defn}

It is sometimes convenient to convert a local function $F$ into
another function with a severely restricted domain and range.

\begin{defn}[Restricted local functions]
If $i$ is a vertex with $k$ neighbors in $Y$, then corresponding to
each function $F$ that is $Y$-local at $i$, we define a function
$f\colon\F^{k+1}\to \F$ where the domain is restricted to the coordinates
corresponding to $i$ and its neighbors, and the output is the new
value $F$ would assign to the $i^{\rm th}$ coordinate under these
conditions.  It should be clear that $F$ and $f$ contain the same
information but packaged in different ways.  Each determines the other
and both have their uses.  Functions such as $F$ can be readily
composed, but functions such as $f$ are easier to describe explicitly
since irrelevant and redundant information has been eliminated.
\end{defn}

The local functions that are easiest to describe are those with extra
symmetries.

\begin{defn}[Symmetric and quasi-symmetric rules]
Let $i$ be a vertex in $Y$ with $k$ neighbors, let $F\colon\F^n\to
\F^n$ be a $Y$-local function at $i$ and let $f\colon\F^{k+1}\to\F$ be its
restricted form. If the output of $f$ only depends on the multiset of
inputs and not their order, in other words, if the states of $i$ and
its neighbors can be arbitrarily permuted without changing the output
of $f$, then $f$ (and $F$) are called \emph{symmetric} local
functions.  If they satisfy the weaker condition that at least the
states of the neighbors of $i$ can be arbitrarily permuted without
changing the output, then $F$ and $f$ are \emph{quasi-symmetric}.  A
list of local functions $\Fy$ is symmetric or quasi-symmetric
when every function in the list has this property.
\end{defn}

The last component of an \sds is an update order.

\begin{defn}[Update orders]
An \emph{update order} $\omega$ is a finite sequence of numbers chosen
from the set $\{1,\ldots, n\}$ such that every number $1 \leq i \leq
n$ occurs at least once.  If every number $1 \leq i \leq n$ occurs
exactly once, then the update order is \emph{simple} (also called a
\emph{permutation}).  Let $W_Y$ denote the collection of all update
orders and let $S_Y$ denote the subset of simple update orders.  The
subscript $Y$ indicates that we are thinking of the numbers in these
sequences as vertices in the graph $Y$.  When considering an arbitrary
update order, we tend to use the notation $\omega = (\omega_1,
\omega_2, \ldots, \omega_m)$ with $m=|\omega|$ (and $m \geq n$, of
course), but when we restrict our attention to simple update orders,
we switch to the notation $\pi = (\pi_1, \pi_2, \ldots, \pi_n)$.
\end{defn}

\begin{defn}[Sequential dynamical systems]
A \emph{sequential dynamical system}, or \sds, is a triple
$(Y,\Fy,\omega)$ consisting of an undirected graph $Y$, a list of
local functions $\Fy$, and an update order $\omega\in W_Y$.  If
$\omega$ is the sequence $(\omega_1,\omega_2,\ldots,\omega_m)$, then
we construct the \emph{\sds map} $[\Fy,\omega]\colon \F^n \to
\F^n$ as the composition $[\Fy,\omega]:=F_{\omega_m} \circ
\cdots \circ F_{\omega_1}$.
\end{defn}

With the usual abuse of notation, we sometimes let the \sds map
$[\Fy,\omega]$ stand in for the entire \sds.  The goal is to
study the behavior of the map $[\Fy,\omega]$ under iteration.  In
this article we focus on the set of states in $\F^n$ that are periodic
and we use $\per[\Fy,\omega] \subset \F^n$ to denote this
collection of periodic states.  The set of periodic states is of
interest both because it is the codomain of high iterates of the \sds
map and the largest subset of states that are permuted by the map.

\begin{defn}[$\omega$- and $\pi$-independence]
A list of $Y$-local functions $\Fy$ is called
\emph{$\omega$-independent} if $\per[\Fy,\omega] =
\per[\Fy,\omega']$ for all update orders $\omega, \omega' \in
W_Y$ and \emph{$\pi$-independent} if $\per[\Fy,\pi] =
\per[\Fy,\pi']$ for all simple update orders $\pi, \pi' \in
S_Y$.
\end{defn}

Every $\omega$-independent $\Fy$ is trivially $\pi$-independent. As we
will see later, the converse does not hold. In the context of this
paper, we will refer to functions that are $\pi$-independent without
being $\omega$-independent as being \emph{exceptional}.
%
%
%
Even though $\pi$-independence is too strong to expect generically,
there are nonetheless many interesting classes of \sdss that have this
property, including two classes where $\pi$-independence is
relatively easy to establish.

\begin{prop}\label{prop:fixed}
If for every simple update order $\pi \in S_Y$, every state in
$\per[\Fy,\pi]$ is fixed by the \sds map $[\Fy,\pi]$,
then $\Fy$ is $\pi$-independent. 
\end{prop}

\begin{proof}
If $\y$ is fixed by $[\Fy,\pi]$, then $\y$ is fixed by each
$F_i$ in $\Fy$ (the simplicity of $\pi$ means that were $F_i$
to change the $i^{\rm th}$ coordinate, there would not be an
opportunity for it to change back).  Being fixed by each $F_i$, $\y$
is also fixed by $[\Fy,\omega]$ for all $\omega \in W_Y$, which
includes all of $S_Y$.  Since this argument is reversible, the \sds
maps with simple update orders share a common set of fixed states.
If, as hypothesized, these are the only periodic states for these
maps, then $\Fy$ is $\pi$-independent.
\end{proof}

In our second example, $\pi$-independence is essentially immediate.

\begin{prop}[Bijective functions]\label{prop:bijective}
If every local function $F_i$ in $\Fy$ is a bijection, then for
every update order $\omega \in W_Y$, $\per[\Fy,\omega]=\F^n$.
As a consequence $\Fy$ is $\omega$-independent, and hence $\pi$-independent.
\end{prop}

\begin{proof}
Since every $F_i$ is a bijection, so is the \sds map $[\Fy,\omega]$
and a sufficiently high iterate is the identity permutation.
\end{proof}

This last result highlights the fact that $\omega$- and $\pi$-independence
focuses on sets rather than cycles, since $\pi$-independent \sdss
with different update orders quite often organize their common
periodic states into different cycle configurations.  In fact, the
restrictions of $\pi$-independent \sds maps with different update
orders to their common periodic states can be used to construct a
group encoding the possible dynamics over this set~\cite{Hansson:05b}.

Collections of $\pi$-independent \sdss also form a natural starting
point for the study of stochastic sequential dynamical
systems. Stochastic finite dynamical systems are often studied through
Markov chains over their state space but in general this leads to
Markov chains with exponentially many states as measured by the number
of cells or vertices.  For $\pi$-independent \sdss one is typically
able to reduce the number of states in such a Markov chain
significantly, at least when focusing on their periodic behavior.

\section{Asynchronous Cellular Automata}\label{sec:aca}

Some of the simplest (classical) cellular automata are the
one-dimensional \cas known as \emph{elementary cellular automata}. In
an elementary \ca, every vertex has precisely two neighbors, the only
possible vertex states are $0$ or $1$, and all local functions are
identically defined.  Since every vertex has two neighbors, the
underlying graph is either a line or a circle and the restricted form
of its common local function is a map $f\colon \F^3\to \F$ where $\F=
\F_2 = \{0,1\}$ is the field with two elements.  There are $2^8 = 256$
such functions, known as \emph{Wolfram rules}, and thus $256$ types of
elementary cellular automata.  Even in such a restrictive situation
there are many interesting dynamical effects to be observed.  The
focus here is on the sequential dynamical systems that correspond to
these classical elementary cellular automata.

Let $Y=\Circ_n$ denote a circular graph with $n$ vertices labeled
consecutively from $1$ to $n$, and to avoid trivialities assume $n>3$.
(The sequential nature of the update rules in an \sds makes infinite
graphs such as lines unsuitable in this context.) Since these are the
only graphs considered in the remainder of the article, we replace
notations such as $W_Y$ or $S_Y$ with $W_n$ and $S_n$, etc.  In
$\Circ_n$ we view the vertex labels as residue classes mod $n$ so that
there is an edge connecting $i$ to $i+1$ for every $i$.

\begin{defn}[Wolfram rules]
Let $F_i\colon \F^n\to \F^n$ be a $\Circ_n$-local function at $i$ and
let $f_i\colon \F^3\to \F$ be its restricted form.  Since the
neighbors of $i$ are $i-1$ and $i+1$, it is standard to list these
coordinates in ascending order in $\F^3$.  Thus, a state $\y\in \F^n$
corresponds to a triple $(y_{i-1},y_i,y_{i+1})$ in the domain of
$f_i$.  Call this a \emph{local state configuration} and keep in mind
that all subscripts are viewed mod $n$.  In order to completely
specify the function $F_i$ it is sufficient to list how the $i^{\rm
  th}$ coordinate is updated for each of the $8$ possible local state
configurations.  More specifically, let $(y_{i-1},y_i,y_{i+1})$ denote
a local state configuration and let $(y_{i-1},z_i,y_{i+1})$ be the
local state configuration after applying $F_i$.  The local function
$F_i$, henceforth referred to as a \emph{Wolfram rule}, is completely
described by the following table.
\begin{equation}
  \label{eqn:wolframdef}
  \begin{array}{c||c|c|c|c|c|c|c|c}
    y_{i-1}y_iy_{i+1} & 111 & 110 & 101 & 100 & 011 & 010 & 001 & 000
    \\ \hline z_i & a_7 & a_6 & a_5 & a_4 & a_3 & a_2 & a_1 & a_0
  \end{array}
\end{equation}
More concisely, the $2^8=256$ possible Wolfram rules can be indexed by
an $8$-digit binary number $a_7a_6a_5a_4a_3a_2a_1a_0$, or by its
decimal equivalent $k = \sum_{i=0}^7 a_i 2^i$.  There is thus one
\emph{Wolfram rule $k$} for each integer $0\leq k \leq 255$. For each
such $n$, $k$ and $i$ let $\Wolf^{(k)}_i$ denote the $\Circ_n$-local
function $F_i\colon \F^n\to \F^n$ just defined, let $\wolf^{(k)}_i$
denote its restricted form $f_i\colon \F^3\to\F$, and let
$\WOLF^{(k)}_n$ denote the list of local functions $(\Wolf^{(k)}_1,
\Wolf^{(k)}_2, \ldots, \Wolf^{(k)}_n)$.  We say that Wolfram rule~$k$
is $\omega$-independent ($\pi$-independent) whenever $\WOLF_n^{(k)}$ is
$\omega$-independent ($\pi$-independent) for all $n>3$.
\end{defn}

For each update order $\omega$ there is an \sds
$(\Circ_n,\WOLF^{(k)}_n,\omega)$ that can be thought of as an
elementary \ca, but with the update functions applied asynchronously
(and possibly more than once).  For this reason, such systems are
called \emph{asynchronous cellular automata} or \acas.  We now state
our main result.

\begin{thm}\label{thm:104}
There are exactly $104$ Wolfram rules that are $\pi$-independent, $86$
of which are additionally $\omega$-independent.  More precisely,
$\WOLF_n^{(k)}$ is $\pi$-independent for all $n>3$ iff $k \in \{$0, 1,
4, 5, 8, 9, 12, 13, \textbf{28}, \textbf{29}, \textbf{32},
\textbf{40}, 51, 54, 57, 60, 64, 65, 68, 69, \textbf{70}, \textbf{71},
72, 73, 76, 77, 78, 79, 92, 93, 94, 95, \textbf{96}, 99, 102, 105,
108, 109, 110, 111, 124, 125, 126, 127, 128, 129, 132, 133, 136, 137,
140, 141, 147, 150, \textbf{152}, 153, 156, \textbf{157}, 160, 164,
168, 172, \textbf{184}, \textbf{188}, 192, 193, \textbf{194}, 195,
196, 197, 198, \textbf{199}, 200, 201, 202, 204, 205, 206, 207, 216,
218, 220, 221, 222, 223, 224, \textbf{226}, 228, \textbf{230}, 232,
234, \textbf{235}, 236, 237, 238, 239, 248, \textbf{249}, 250,
\textbf{251}, 252, 253, 254, 255$\}$. The rules in boldface are the
$18$ exceptional rules -- they are $\pi$-independent but not
$\omega$-independent.
\end{thm}

The main result of \cite{Hansson:05b} states that precisely 11 of the
16 \emph{symmetric} Wolfram rules are $\omega$-independent over
$\Circ_n$ for all $n>3$. Theorem~\ref{thm:104} significantly extends
this result, reproving it in the process.  In addition to identifying
a large class of $\omega$- and $\pi$-independent \acas, the proof also
provides further insight into the dynamics of these systems at both
periodic and transient states and thus serves as a foundation for the
future study of their stochastic properties.  We conclude this section
with two remarks about the role played by computer investigations of
these systems.

\begin{rem}[Unlisted numbers]
The ``only if'' portion of this theorem was established
experimentally.  For each $4\leq n \leq 9$, for each $0 \leq k \leq
255$, and for each simple update order $\pi \in S_n$, a computer
program written by the first and third authors calculated the set
$\per[\WOLF_n^{(k)},\pi]$.  For each of the 152 values of $k$ not
listed above, there were distinct simple update orders that led to
distinct sets of periodic states, leaving the remaining 104 rules as
the only ones with the potential to be $\pi$-independent for all
$n>3$.  Moreover, since a counterexample for one value of $n$ leads to
similar counterexamples for all multiples of $n$, these $104$ rules
are also the only ones that are eventually $\pi$-independent for all
sufficiently large values of $n$.  Because these brute-force
calculations are explicit yet tedious they have been omitted, but the
interested reader should feel free to contact the third author for a
copy of the software that performed the calculations.
\end{rem}

\begin{rem}[Computational guidance]
These early computer-aided investigations also had a major impact on
the ``if'' portion of the proof.  Once the computer results
highlighted the 104 rules that were $\pi$-independent for small
values of $n$, we identified patterns and clusters among the 104
rules, which led to conjectured lemmas, and eventually to proofs that
our conjectures were correct.  The computer calculations thus provided
crucial data that both prompted ideas and tempered our search for
intermediate results.
\end{rem}

\section{Wolfram Rule Notations}\label{sec:notation}

Patterns among the 104 numbers listed in Theorem~\ref{thm:104} are
difficult to discern because the conversion from binary to decimal
obscures many structural details.  In this section we introduce other
ways to describe the Wolfram rules that makes their similarities and
differences more immediately apparent.

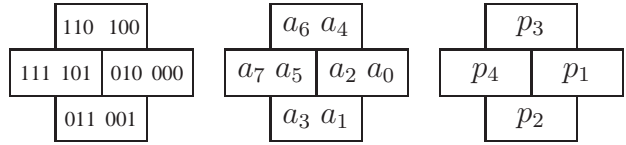
\begin{figure}[t]
\begin{tabular}{ccc}
  \setlength{\unitlength}{1.2cm}
  \begin{picture}(2,1.7)
    \put(.5,1.015){\framebox(1,.49){{\tiny 110} \boxsp\!\!\!\! {\tiny 100}}}
    \put(0, .5){\framebox(1,.5){{\tiny 111} \boxsp\!\!\!\!\! {\tiny 101}}}
    \put(1.01, .5){\framebox(1,.5){{\tiny 010} \boxsp\!\!\!\!\! {\tiny 000}}}
    \put(.5,0){\framebox(1,.485){{\tiny 011} \boxsp\!\!\!\!\! {\tiny 001}}}
  \end{picture}
&
  \setlength{\unitlength}{1.2cm}
  \begin{picture}(2,1.7)
    \put(.5,1.015){\framebox(1,.49){$a_6$ \boxsp\!\!\!\!\! $a_4$}}
    \put(0, .5){\framebox(1,.5){$a_7$ \boxsp\!\!\!\!\! $a_5$}}
    \put(1.01, .5){\framebox(1,.5){$a_2$ \boxsp\!\!\!\!\! $a_0$}}
    \put(.5,0){\framebox(1,.485){$a_3$ \boxsp\!\!\!\!\! $a_1$}}
  \end{picture}
&
\setlength{\unitlength}{1.2cm}
    \begin{picture}(2,1.7)
      \put(.5,1.015){\framebox(1,.49){$p_3$}}
      \put(0, .5){\framebox(1,.5){$p_4$}}
      \put(1.01, .5){\framebox(1,.5){$p_1$}}
      \put(.5,0){\framebox(1,.485){$p_2$}}
    \end{picture}
\end{tabular}
\caption{Grid notation for Wolfram rules\label{fig:grid}}
\end{figure}

\begin{figure}[t]
\begin{tabular}{cc}
  \setlength{\unitlength}{1.2cm}
  \begin{picture}(2,1.7)
    \put(.5,1.015){\framebox(1,.49){0 \boxsp 1}}
    \put(0, .5){\framebox(1,.5){0 \boxsp 0}}
    \put(1.01, .5){\framebox(1,.5){1 \boxsp 1}}
    \put(.5,0){\framebox(1,.485){1 \boxsp 0}}
  \end{picture}
&
\setlength{\unitlength}{1.2cm}
    \begin{picture}(2,1.7)
      \put(.5,1.015){\framebox(1,.49){{\tt x}}}
      \put(0, .5){\framebox(1,.5){{\tt 0}}}
      \put(1.01, .5){\framebox(1,.5){{\tt 1}}}
      \put(.5,0){\framebox(1,.485){{\tt -}}}
    \end{picture}
\end{tabular}
\caption{Converting Wolfram rule 29 = 00011101 into its tag {\tt
    0x-1}.\label{fig:tag-conversion}}
\end{figure}
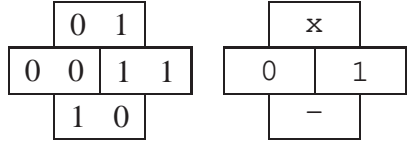

\begin{defn}[Grid notation]
For each binary number $k=a_7a_6a_5a_4a_3a_2a_1a_0$ we arrange its
digits in a grid.  The $8$ local state configurations can be viewed as
the vertices of a $3$-cube and we arrange them according to the
conventional projection of a $3$-cube into the plane.  See the
left-hand side of Figure~\ref{fig:grid}.  Next, we can place the binary
digits of $k$ at these positions as shown in the center of
Figure~\ref{fig:grid}.  The boxes have been added because the local
state configurations come in pairs.  When a local function is applied,
the states of the neighbors of $i$ are left unchanged, so that the
resulting local state configuration is located in the same box.  We
call this the \emph{grid notation} for $k$.  The grid notation for
Wolfram rule 29 = 00011101 is shown on the left-hand side of
Figure~\ref{fig:tag-conversion}.
\end{defn}

Because grid notation is sometimes cumbersome to work with we also
define a very concise $4$ symbol tag for each Wolfram rule that
respects the box structure of the grid.

\begin{defn}[Tags]
When we look at the grid notation for a Wolfram rule, in each box we
see a pair of numbers, $11$, $00$, $10$, or $01$, and we encode these
configurations by the symbols {\tt 1}, {\tt 0}, {\tt -}, and {\tt x},
respectively.  In other words
`{\tt 1}' = \framebox[1cm]{$1$ \boxsp $1$}\ ,
`{\tt 0}' = \framebox[1cm]{$0$ \boxsp $0$}\ ,
`{\tt -}' = \framebox[1cm]{$1$ \boxsp $0$}\ , and
`{\tt x}' = \framebox[1cm]{$0$ \boxsp $1$}.
The symbols are meant to indicate that when the states of the
neighbors place us in this box, the local function updates the $i^{\rm
  th}$ coordinate by converting it to a $1$, converting it to a $0$,
leaving it unchanged, or always changing it.  We label the symbols for
the four boxes $p_1$, $p_2$, $p_3$ and $p_4$ as shown on the right-hand
side of Figure~\ref{fig:grid} and we define the \emph{tag} of $k$ to
be the string $p_4p_3p_2p_1$.  The numbering and the order of the
$p_i$'s have been chosen to match the binary representation as closely
as possible, with the hope of easing conversions between binary and
tag representations.  The process of converting Wolfram rule $29$ to
its tag {\tt 0x-1} is illustrated in Figure~\ref{fig:tag-conversion}.
\end{defn}

\begin{defn}[Symmetric and asymmetric]
The middle row of the grid contains the positions where the states of
the neighbors are equal and the top and bottom rows contain the
positions where the states of the neighbors are different.  We call
the middle row the \emph{symmetric} portion of the grid and the top
and bottom rows the \emph{asymmetric} portion.  In the tag
representation, the beginning and end of a tag describes how the rule
responds to a symmetric neighborhood configuration and the middle of a
tag describes how it responds to an asymmetric neighborhood
configuration.  With this in mind we call $p_4p_1$ the \emph{symmetric
  part} of the tag $k = p_4p_3p_2p_1$ and we call $p_3p_2$ its
\emph{asymmetric part}.
\end{defn}

Table~\ref{table:16x10} shows the $104$ $\pi$-independent Wolfram
rules listed in Theorem~\ref{thm:104} arranged according to the
symmetric and asymmetric parts of their tags.  The rows list all 16
possibilities for the symmetric part of the tag while the columns list
only 10 of the 16 possibilities for the asymmetric part since only
these 10 occur among the 104 rules.  In addition, each row and column
label has a decimal equivalent, listed next to the row and column
headings, that add up to $k$.  In this format the benefits of the tag
representation should be clear.  Far from being distributed
haphazardly, the $\pi$-independent rules appear clustered together
in large blocks.
Table~\ref{table:16x10} reveals a lot of structure, but some patterns
remain slightly hidden due to the order in which the rows and columns
are listed.  For example, there is a $4$-by-$4$ block of bijective
rules obtained by restricting attention to the four rows that show up
in the last column and the four columns that show up in the last row.

\begin{table}[t]
\begin{tabular}{c|c||cccccccccc|}
& $p_3$ &{\tt -}&{\tt -}&{\tt 0}&{\tt 0}&{\tt -}&{\tt 1}&{\tt 1}&{\tt
    -}&{\tt x}&{\tt x}\vspace{-.15cm} \\
& $p_2$ &{\tt -}&{\tt 0}&{\tt -}&{\tt 0}&{\tt 1}&{\tt -}&{\tt 1}&{\tt
    x}&{\tt -}&{\tt x}\\
\hline
$p_4p_1$ &  &  72 &  64 &   8 &   0 &  74 &  88 &  90 &  66 &  24 & 18 \\
\hline\hline
{\tt --}&132& 204 & 196 & 140 & 132 & 206 & 220 & 222 & 198 & 156 & 150 \\
{\tt 0-}&  4&  76 &  68 &  12 &   4 &  78 &  92 &  94 &  \textbf{70} &  \textbf{28} &  \\
{\tt -0}&128& 200 & 192 & 136 & 128 & 202 & 216 & 218 & \textbf{194} & \textbf{152} &  \\
{\tt 1-}&164& 236 & 228 & 172 & 164 & 238 & 252 & 254 & \textbf{230} & \textbf{188} &  \\
{\tt -1}&133& 205 & 197 & 141 & 133 & 207 & 221 & 223 & \textbf{199} & \textbf{157} &  \\
{\tt 10}&160& 232 & 224 & 168 & 160 & 234 & 248 & 250 & \textbf{226} & \textbf{184} &  \\
{\tt 01}&  5&  77 &  69 &  13 &   5 &  79 &  93 &  95 &  \textbf{71} &  \textbf{29} &  \\
{\tt 00}&  0&  72 &  64 &   8 &   0 &     &     &     &     &     &  \\
{\tt x0}& 32&     &  \textbf{96} &  \textbf{40} &  \textbf{32} &     &     &     &     &     &  \\
{\tt 0x}&  1&  73 &  65 &   9 &   1 &     &     &     &     &     &  \\
{\tt -x}&129& 201 & 193 & 137 & 129 &     &     &     & 195 & 153 & 147 \\
{\tt x-}& 36& 108 &     &     &     & 110 & 124 & 126 & 102 &  60 & 54 \\
{\tt x1}& 37& 109 &     &     &     & 111 & 125 & 127 &     &     &  \\
{\tt 1x}&161&     &     &     &     & \textbf{235} & \textbf{249} & \textbf{251} &     &     &  \\
{\tt 11}&165& 237 &     &     &     & 239 & 253 & 255 &     &     &  \\
{\tt xx}& 33& 105 &     &     &     &     &     &     &  99 &  57 & 51 \\
\hline
\end{tabular}
\vspace*{3mm}
\caption{The $104$ $\pi$-independent Wolfram rules arranged by the
  symmetric and asymmetric parts of their tags. As before, the $18$
  exceptional rules are in boldface.
  \label{table:16x10}}
\end{table}

\begin{prop}[Bijective rules]\label{prop:bijective-rules}
Wolfram rules 51, 54, 57, 60, 99, 102, 105, 108, 147, 150, 153, 156,
195, 198, 201 and 204 are $\omega$-independent.
\end{prop}

\begin{proof}
The 16 rules listed have tags where each $p_i$ is either {\tt -} or
{\tt x}.  These (and only these) Wolfram rules correspond to bijective
local functions and by Proposition~\ref{prop:bijective} the \acas
these rules define are $\omega$-independent.
\end{proof}

\section{Dynamical equivalence}\label{sec:equivalence}

In this section we use the notion of dynamical equivalence to reduce
the proof of Theorem~\ref{thm:104} to a more manageable size.  Two
sequential dynamical systems $(Y,\Fy,\omega)$ and
$(Y,\frak{F}'_Y,\omega')$ defined over the same graph $Y$ are said to
be \emph{dynamically equivalent} if there is a bijection
$H\colon\F^n\to\F^n$ between their states such that $H \circ [
\Fy,\omega ] = [\frak{F}'_Y,\omega' ] \circ H$.  The key fact about
dynamically equivalent \sdss, which is also easy to show, is that $H$
establishes a bijection between their periodic states.  In particular,
$H(\per [\Fy,\omega]) = \per[\frak{F}'_Y,\omega']$.  Thus, if
$\frak{F}'_Y$ is an $\omega$-independent \sds and for each $\omega\in
W_Y$ there exists an $\omega'\in W_Y$ such that $(Y,\Fy,\omega)$ and
$(Y,\frak{F}'_Y,\omega')$ are dynamically equivalent using \emph{the
same function} $H$, then $\Fy$ is also $\omega$-independent. An
analogous statement about $\pi$-independence holds -- restict to the
simple update orders and replace $\omega$ with $\pi$ above.

Although there are $256$ Wolfram rules, many give rise to dynamically
equivalent \acas.  In particular, there are three relatively
elementary ways to alter an \aca to produce another one that appears
different on the surface, but which is easily seen to be dynamically
equivalent to the original.  These are obtained by (1) renumbering the
vertices in the opposite direction, (2) systematically switching all
$1$s to $0$s and $0$s to $1$s, or (3) doing both at once.  We call
these alterations \emph{reflection}, \emph{inversion} and
\emph{reflection-inversion} of the \aca, respectively.  The term
reflection highlights the fact that this alteration makes it appear as
though we picked up the circular graph and flipped it over.  We begin
by describing the effect renumbering has on individual local
functions.

\begin{defn}[Renumbering]
The renumbering of the vertices we have in mind is achieved by the map
$r\colon \Circ_n \to \Circ_n$ that sends vertex $i$ to vertex $n+1-i$.
For later use we extend this to a map $r\colon W_n\to W_n$ on update
orders by applying $r$ to each entry in the sequence.  More
specifically, if $\omega = (\omega_1, \omega_2, \ldots, \omega_m)$,
then $r(\omega) = (r(\omega_1), r(\omega_2), \ldots, r(\omega_m))$.
Finally, on the level of states we define a map $\R\colon \F^n\to
\F^n$ that sends sends $\y = (y_1, y_2, \ldots, y_n)$ to $(y_n,
\ldots, y_2, y_1)$, and we note that $\R$ is an involution.
\end{defn}

\begin{defn}[Reflected rules]
If the vertices of $\Circ_n$ are renumbered, rule $\Wolf_i^{(k)}$ is
applied, and then the renumbering is reversed, the net effect is the
same as if a different Wolfram rule were applied to the vertex $r(i)$.
Let $\ell$ be the number that represents this other Wolfram rule.  The
differences between $k$ and $\ell$ are best seen in grid notation.
The renumbering not only changes the vertex at which the rule seems to
be applied, but it also reverses the order in which the $3$
coordinates are listed in the restricted local form.  Only the
asymmetric local state configurations, i.e. the top and bottom rows of
the grid, are altered by this change so that the grid for $\ell$ looks
like a reflection of the grid for $k$ across a horizontal line.  We
call $\ell$ the \emph{reflection} of $k$ and we define a map
$\refl\colon \{0,\ldots,255\}\to \{0,\ldots,255\}$ with $\refl(k) =
\ell$.  On the level of tags, the only change is to switch order of
$p_2$ and $p_3$, so, for example $\ell$={\tt 01-x} is the reflection
of $k$={\tt0-1x}.
\end{defn}

In short, when $\ell = \refl(k)$, $\R \circ \Wolf_i^{(k)} \circ \R =
\Wolf_{r(i)}^{(\ell)}$ and, since $\R$ is an involution, this can be
rewritten as $\R \circ \Wolf_i^{(k)} = \Wolf_{r(i)}^{(\ell)} \circ
\R$.

\begin{prop}\label{prop:refl}
If $\ell = \refl(k)$, then $\WOLF_n^{(k)}$ is $\omega$-independent
($\pi$-independent) iff $\WOLF_n^{(\ell)}$ is $\omega$-independent
($\pi$-independent).
\end{prop}

\begin{proof}
The value of $\ell$ was defined so that $\R \circ \Wolf_i^{(k)} =
\Wolf_{r(i)}^{(\ell)} \circ \R$.  As a result, for any $\omega \in
W_n$, the \aca $(\Circ_n,\WOLF_n^{(k)},\omega)$ is dynamically
equivalent to the \aca $(\Circ_n,\WOLF_n^{(\ell)},r(\omega))$ since
\begin{equation*}
\begin{array}{rcl}
\R \circ [\WOLF_n^{(k)},\omega] & = & \R \circ \Wolf_{\omega_m}^{(k)}
\circ \cdots \circ \Wolf_{\omega_2}^{(k)} \circ
\Wolf_{\omega_1}^{(k)}\\ &=& \Wolf_{r(\omega_m)}^{(\ell)} \circ \cdots
\circ \Wolf_{r(\omega_2)}^{(\ell)} \circ \Wolf_{r(\omega_1)}^{(\ell)}
\circ \R\\ &=& [\WOLF_n^{(\ell)},r(\omega)] \circ \R.
\end{array}
\end{equation*}
The argument at the beginning of the section now shows that the
$\omega$-independence of $\WOLF_n^{(\ell)}$ implies that of
$\WOLF_n^{(k)}$, but since $\ell=\refl(k)$ implies $k=\refl(\ell)$,
the converse also holds. The argument for $\pi$-independence is
analogous.
\end{proof}

Similar results hold for inversions as we now show.

\begin{defn}[Inverting]
Let $\one$ and $\zer$ denote the special states $(1,1,\ldots,1)$ and
$(0,0,\ldots,0)$ in $\F^n$.  Since the function $i(a)=1-a$ changes $1$
to $0$ and $0$ to $1$, the map $\I\colon \F^n\to \F^n$ sending $\y$ to
$\one - \y$, has this effect on each coordinate of $\y$. The map $\I$
is an involution like $\R$, and from their definitions it is easy to
check that they commute with each other.
\end{defn}

\begin{defn}[Inverted rules]
If the states of $\Circ_n$ are inverted, rule $\Wolf_i^{(k)}$ is
applied, and then the inversion is reversed, the net effect is the
same as if a different Wolfram rule were applied at vertex $i$.  Let
$\ell$ be the number that represents this other Wolfram rule.  The
differences between $k$ and $\ell$ are again best seen in grid
notation.  The pre-inversion of states effects the local state
configurations as though the grid had been rotated $180^\circ$.  The
second inversion merely changes every entry so that $1$s become $0$s
and $0$s become $1$s.  Thus the grid for $\ell$ can be obtained from
the grid for $k$ by rotating the grid and altering every entry.  We
call $\ell$ the \emph{inversion} of $k$ and define a map
$\inv\colon\{0,\ldots,255\}\to \{0,\ldots,255\}$ with $\inv(k) =
\ell$.  On the level of tags, there are two changes that take place.
Boxes $p_1$ and $p_4$ switch places as do boxes $p_2$ and $p_3$, but
in process the boxes are turned over and the numbers changed.  If we
look at what this does to the entries in a box, $11$ becomes $00$,
$00$ becomes $11$, while $10$ and $01$ are left unchanged.  To
formalize this, define a conjugation map $c\colon \{{\tt 1,0,-,x}\}
\to \{{\tt 1,0,-,x}\}$ with $c({\tt 1}) ={\tt 0}$, $c({\tt 0}) ={\tt
  1}$, $c({\tt -}) ={\tt -}$, and $c({\tt x}) ={\tt x}$.  When $k$ has
tag $p_4p_3p_2p_1$, $\ell$ has tag $c(p_1) c(p_2) c(p_3) c(p_4)$, so,
for example, $\ell$ = {\tt x0-1} is the inversion of $k$ = {\tt 0-1x}.
\end{defn}

In short, when $\ell = \inv(k)$, $\I \circ \Wolf_i^{(k)} \circ \I =
\Wolf_i^{(\ell)}$ and, since $\I$ is an involution, this can be
rewritten as $\I \circ \Wolf_i^{(k)} = \Wolf_i^{(\ell)} \circ \I$.

\begin{prop}\label{prop:inv}
If $\ell = \inv(k)$, then $\WOLF_n^{(k)}$ is $\omega$-independent
($\pi$-independent) iff $\WOLF_n^{(\ell)}$ is $\omega$-independent
($\pi$-independent).
\end{prop}

\begin{proof}
The value of $\ell$ was defined so that $\I \circ \Wolf_i^{(k)} =
\Wolf_i^{(\ell)} \circ \I$.  As in the proof of
Proposition~\ref{prop:refl} this implies that for any $\omega \in
W_n$, the \aca $(\Circ_n,\WOLF_n^{(k)},\omega)$ is dynamically
equivalent to the \aca $(\Circ_n,\WOLF_n^{(\ell)},\omega)$.  The
argument at the beginning of the section and the fact that
$\ell=\inv(k)$ implies $k=\inv(\ell)$, complete the proof as
before. The argument for $\pi$-independence is completely analogous.
\end{proof}

As an immediate corollary of Propositions~\ref{prop:refl}
and~\ref{prop:inv}, when $\ell = \refl(\inv(k)) = \inv(\refl(k))$,
$\WOLF_n^{(k)}$ is $\omega$-independent ($\pi$-independent) iff
$\WOLF_n^{(\ell)}$ is $\omega$-independent ($\pi$-independent). If we
partition the $256$ Wolfram rules into equivalence classes of rules
related by reflection, inversion or both, then there are $88$ distinct
equivalence classes and the $104$ rules listed in
Theorem~\ref{thm:104} are the union of $41$ of them.

\begin{figure}[t]
\begin{tabular}{ccc}
\begin{tabular}{c|c||cc|}
& $p_3$ &{\tt 0}&{\tt 0}\\
& $p_2$ &{\tt 0}&{\tt -}\\
\hline
$p_4p_1$ &  & 0 & 8 \\
\hline\hline
{\tt --} & 132 & 132 & 140 \\
{\tt 0-} &   4 &   4 &  12 \\
{\tt -0} & 128 & 128 & 136 \\
{\tt 00} &   0 &   0 &   8 \\
\hline
{\tt -1} & 133 & 133 & 141 \\
{\tt 01} &   5 &   5 &  13 \\
{\tt -x} & 129 & 129 & 137 \\
{\tt 0x} &   1 &   1 &   9 \\
\hline
{\tt 1-} & 164 & 164 & 172 \\
{\tt 10} & 160 & 160 & 168 \\
{\tt x0} &  32 &  \textbf{32} &  \textbf{40} \\
\hline
\end{tabular}
& \hspace*{1cm} &
\begin{tabular}{c|c||ccc|}
& $p_3$ &{\tt -}&{\tt x}&{\tt x}\\
& $p_2$ &{\tt -}&{\tt -}&{\tt x}\\
\hline
$p_4p_1$ &  &  72 &  24 & 18 \\
\hline\hline
{\tt --}&132& 204 & 156 & 150\\
{\tt x-}& 36& 108 &  60 & 54 \\
{\tt xx}& 33& 105 &  57 & 51\\
\hline
{\tt -0}&128& 200 & \textbf{152} &  \\
{\tt 10}&160& 232 & \textbf{184} &  \\
\hline
{\tt 0-}&  4&  76 &  \textbf{28} &  \\
{\tt 01}&  5&  77 &  \textbf{29} &  \\
{\tt 00}&  0&  72 &     &  \\
{\tt 0x}&  1&  73 &     &  \\
\hline
\end{tabular}
\vspace*{5mm}
\end{tabular}
\caption{The $41$ $\pi$-independent Wolfram rules up to
  equivalence, separated into two tables by their behavior in
  asymmetric contexts. The exceptional rules are shown in bold.\label{fig:41}}
\end{figure}

Figure~\ref{fig:41} displays representatives of these $41$ classes in
pared down versions of Table~\ref{table:16x10}.  We used reflection
and inversion to eliminate $5$ of the $10$ columns.  Every rule with a
{\tt 1} in the asymmetric portion of its tag is the inversion of a rule
with a {\tt 0} instead.  In particular, the entries in the $3$ columns
headed {\tt -1}, {\tt 1-} and {\tt 11} are inversions of the entries in
the columns headed {\tt 0-}, {\tt -0} and {\tt 00}, respectively.
Next, since reflections switch $p_2$ and $p_3$ we can also eliminate
the columns headed {\tt -0}, {\tt -x} as redundant.  This leaves the
$5$ columns headed {\tt 00}, {\tt 0-}, {\tt --}, {\tt x-} and {\tt
  xx}.  Since the last $3$ do not contain {\tt 0}s or {\tt 1}s,
further inversions, or inversion-reflections can be used to identify
redundant rows in these columns.

As mentioned above, the $41$ rules listed in Figure~\ref{fig:41} are
representatives of the $41$ distinct equivalence classes of rules
whose $\pi$-independence needs to be established in order to prove
Theorem~\ref{thm:104}.  The rows in each table have been arranged to
correspond as closely as possible with the structure of the proof.
For example, the first three rows of the table on the right-hand side
are the $9$ equivalence classes shown to be $\pi$-independent by
Proposition~\ref{prop:bijective}.

\section{Major Classes}\label{sec:major}

In this section we prove that four large sets of Wolfram rules are
$\omega$-independent.  All of the proofs are similar and, when
combined with Proposition~\ref{prop:bijective}, leave only the $6$
equivalence classes of the exceptional Wolfram rules that need to be
discussed separately. The main tool we use in this section is the
notion of a potential function.

\begin{defn}[Potential functions]
Let $F\colon X \to X$ be a map whose dynamics we wish to understand.
A \emph{potential function for $F$} is any map $\rho\colon X\to \real$ such
that $\rho(F(x))\leq \rho(x)$ for all $x\in X$.  A potential function
narrows our search for periodic points since any element $x$ with
$\rho(F(x)) < \rho(x)$ cannot be periodic: further applications of $F$
can never return $x$ to its original potential, hence the name.  The
only elements in $X$ that are possibly periodic under $F$ are those
whose potential under $\rho$ never drops at all.  If we call the
inverse image of a number in $\real$ a \emph{level set of $\rho$},
then to find all periodic points of $F$, we only need to examine its
behavior on each of these level sets. Finally, it should be clear that
when non-decreasing functions are used in the definition instead of
non-increasing ones, the effect is the same.
\end{defn}

\begin{defn}[\sds potential functions]
A potential function for an \sds such as $(Y,\Fy,\omega)$ is a
map $\rho\colon\F^n\to \real$ that is a potential function, in the sense
defined above, for the \sds map $[\Fy,\omega]$. The easiest way
to create such a function is to find one that is a potential function
for every local function $F_i$ in $\Fy$.  Of course, $\rho$
should be either a non-decreasing potential function for each $F_i$ or
a non-increasing potential function for each $F_i$, rather than a
mixture of the two, for the inequalities to work out.  When $\rho$ has
this stronger property we call it a \emph{potential function for}
$\Fy$ since such a $\rho$ is a potential function for
$(Y,\Fy,\omega)$ for every choice of update order $\omega$.
\end{defn}

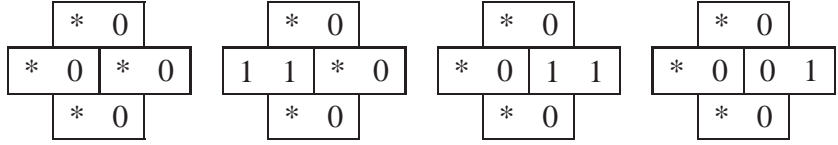
\begin{figure}[t]
\begin{tabular}{cccc}
\setlength{\unitlength}{1.2cm}
    \begin{picture}(2,1.7)
      \put(.5,1.015){\framebox(1,.49){* \boxsp 0}}
      \put(0, .5){\framebox(1,.5){* \boxsp 0}}
      \put(1.01, .5){\framebox(1,.5){* \boxsp 0}}
      \put(.5,0){\framebox(1,.485){* \boxsp 0}}
    \end{picture}
& 
\setlength{\unitlength}{1.2cm}
    \begin{picture}(2,1.7)
      \put(.5,1.015){\framebox(1,.49){* \boxsp 0}}
      \put(0, .5){\framebox(1,.5){1 \boxsp 1}}
      \put(1.01, .5){\framebox(1,.5){* \boxsp 0}}
      \put(.5,0){\framebox(1,.485){* \boxsp 0}}
    \end{picture}
& 
\setlength{\unitlength}{1.2cm}
    \begin{picture}(2,1.7)
      \put(.5,1.015){\framebox(1,.49){* \boxsp 0}}
      \put(0, .5){\framebox(1,.5){* \boxsp 0}}
      \put(1.01, .5){\framebox(1,.5){1 \boxsp 1}}
      \put(.5,0){\framebox(1,.485){* \boxsp 0}}
    \end{picture}
& 
\setlength{\unitlength}{1.2cm}
    \begin{picture}(2,1.7)
      \put(.5,1.015){\framebox(1,.49){* \boxsp 0}}
      \put(0, .5){\framebox(1,.5){* \boxsp 0}}
      \put(1.01, .5){\framebox(1,.5){0 \boxsp 1}}
      \put(.5,0){\framebox(1,.485){* \boxsp 0}}
    \end{picture}
\end{tabular}
\caption{Four major classes of $\omega$-independent Wolfram
  rules.\label{fig:major-forms}}
\end{figure}

\begin{prop}\label{prop:one-way}
Rules 0, 4, 8, 12, 72, 76, 128, 132, 136, 140 and 200 are
$\omega$-independent.
\end{prop}

\begin{proof}
If $k$ is one of the numbers listed above, then its grid notation
matches the leftmost form shown in Figure~\ref{fig:major-forms}.
(Each $*$ is to be interpreted as either a $0$ or a $1$ so that 16
rules share this form, the 11 listed in the statement and 5 that are
equivalent to the listed rules or to previously known cases.)  The $4$
specified values mean that local functions never remove $0$s.  Thus,
the map $\rho$ sending $\y\in \F^n$ to the number of $0$s it contains
is a non-decreasing potential function for $\WOLF_n^{(k)}$.  Moreover,
the local functions $\Wolf_i^{(k)}$ cannot change $\y$ without raising
$\rho(\y)$, so all periodic states are fixed states (for any update
order), and by Proposition~\ref{prop:fixed} $\WOLF_n^{(k)}$ is
$\omega$-independent.
\end{proof}

For the next potential function, additional definitions are needed.

\begin{defn}[Blocks]
A state $\y\in\F^n$ is thought of as a cyclic binary $n$-bit string
with the indices taken mod $n$, and a \emph{substring} of $\y$
corresponds to a set of consecutive indices. We refer to maximal
substrings of all $0$s as \emph{0-blocks} and maximal substrings of
all $1$s as \emph{1-blocks}.  If a block contains only a single number
it is \emph{isolated} and if it contains more than one number it is
\emph{non-isolated}.  The state $\y = 010110$, for example contains
one isolated $0$-block and one non-isolated $0$-block of length~$2$
that wraps across the end of the word.
\end{defn}

We study how these blocks evolve as the local functions are applied.
The decomposition of a Wolfram rule into its symmetric and asymmetric
parts is particularly well adapted to the study of these evolutions.
The asymmetric rules either make no change or shrink a non-isolated
$1$-block or $0$-block from the left or the right, depending on which
of the $4$ asymmetric rules we are considering.  Similarly, the $4$
symmetric rules either do nothing, they remove an isolated block or
they create an isolated block in the interior of a long block.

\begin{prop}\label{prop:B}
Rules 160, 164, 168, 172 and 232 are $\omega$-independent.
\end{prop}

\begin{proof}
If $k$ is one of the numbers listed above, then its grid notation
matches the second form shown in Figure~\ref{fig:major-forms}. The
specified values mean that (1) the only $0$s ever removed are the
isolated $0$s and (2) isolated $0$s are never added.  In particular,
non-isolated blocks of $0$s persist indefinitely, they might grow but
they never shrink or split, and the isolated $0$s, once removed, never
return.  Thus, the map $\rho$ that sends $\y$ to the number of
non-isolated $0$s in $\y$ minus the number of isolated $0$s in $\y$ is
a non-decreasing potential function for $\WOLF_n^{(k)}$.  As before,
the local functions $\Wolf_i^{(k)}$ cannot change $\y$ without raising
$\rho(\y)$, so all periodic states are fixed states (for any update
order), and by Proposition~\ref{prop:fixed} $\WOLF_n^{(k)}$ is
$\omega$-independent.
\end{proof}

\begin{prop}\label{prop:C}
Rules 5, 13, 77, 133 and 141 are $\omega$-independent.
\end{prop}

\begin{proof}
If $k$ is one of the numbers listed above, then its grid notation
matches the third form shown in Figure~\ref{fig:major-forms}.  This
time the specified values mean that (1) the only $0$s that are removed
create isolated $1$s, and (2) isolated $1$s are never removed and they
never stop being isolated.  Thus the map $\rho$ that sends $\y$ to the
number of $0$s in $\y$ plus \emph{twice} the number of isolated $1$s
in $\y$ is a non-decreasing potential function for $\WOLF_n^{(k)}$.
Once again, the local functions $\Wolf_i^{(k)}$ cannot change $\y$
without raising $\rho(\y)$, so all periodic states are fixed states
(for any update order), and by Proposition~\ref{prop:fixed}
$\WOLF_n^{(k)}$ is $\omega$-independent.
\end{proof}

The argument for the fourth collection is slightly more complicated.

\begin{prop}\label{prop:D}
Rules 1, 9, 73, 129 and 137 are $\omega$-independent.
\end{prop}

\begin{proof}
If $k$ is one of the numbers listed above, then its grid notation
matches the rightmost form shown in Figure~\ref{fig:major-forms}.
This time the specified values mean that (1) the only $0$s that are
removed create isolated $1$s, but (2) isolated $1$s can also be
removed.  The map $\rho$ that sends $\y$ to the number of $0$s in $\y$
plus the number of isolated $1$s in $\y$ is a non-decreasing potential
function for $\Wolf_n^{(k)}$, but the difficulty is that there are
local changes with $\rho(\Wolf_i^{(k)}(\y)) = \rho(\y)$.  This is true
for the local change $000\to 010$ and for the local change $010\to
000$.  All other local changes raise the potential, but the existence
of these two equalities indicates that there might be (and there are)
states that are periodic under the action of some \sds map
$[\WOLF_n^{(k)},\omega]$ without being fixed. Rather than appeal to a
general theorem, we calculate its periodic states explicitly in this
case.

Fix an update order $\omega\in W_n$ and, for convenience, let
$F\colon\F^n\to \F^n$ denote the \sds map $[\WOLF_n^{(k)},\omega]\colon\F^n\to
\F^n$.  If $a_3=0$ and $\y$ contains a substring of the form $011$,
then $\rho(F(\y)) > \rho(\y)$ and $\y$ is not periodic under $F$.
This is because either (1) the substring remains unaltered until its
central coordinate is updated, at which point it changes to $0$ and
$\rho$ is raised, or (2) it is altered ahead of time by switching the
$1$ on the right to a $0$ (also raising $\rho$), or by switching the
$0$ on the left to a $1$ (impossible since $a_1=a_5=0$).  Analogous
arguments show that if $a_6=0$ and $\y$ contains the substring $110$,
or if $a_7=0$ and $\y$ contains the substring $111$, then $\y$ is not
periodic under $F$.  Let $P$ be the subset of $\F^n$ where these
situations do not occur.  More specifically, if $a_3=0$ remove the
states with $011$ substrings, if $a_6=0$ remove the states with $110$
substrings, and if $a_7=0$ remove the states with $111$ substrings. If
all three are equal to $1$, then $P = \F^n$.

We claim that $P=\per[\WOLF_n^{(k)},\omega]$, independent of the
choice of $\omega$.  We have already shown $P \supset
\per[\WOLF_n^{(k)}, \omega]$.  Note that $P$ is invariant under $F$
(in the sense that $F(P) \subset P$) since the allowed local changes
are not able to create the forbidden substrings when they do not
already exist.  Moreover, $F$ restricted to $P$ agrees with rule 201 =
{\tt ---x}, the rule of this form with $a_3 = a_6 = a_7 = 1$, since
whenever $a_3$, $a_6$ or $a_7$ is $0$, $P$ has been suitably
restricted to make this fact irrelevant.  Finally, for every $\omega$
rule 201 is bijective, thus $F$ is injective on $P$, $F$ permutes the
states in $P$ and a sufficiently high power of $F$ is the identity,
showing every state in $P$ is periodic independent of our choice of
$\omega$.
\end{proof}

\section{Exceptional Cases}\label{sec:minor}

At this point we come to the exceptional rules -- those that are
$\pi$-independent but fail to be $\omega$-independent. They consist of
$6$ rules that come in pairs: 28 and 29, 32 and 40, and 152 and 184.
These final $6$ rules exhibit more intricate dynamics and the proofs
are, of necessity, more delicate. We treat them in order of
difficulty.

\begin{prop}\label{prop:E}
Rules 32 and 40 are $\pi$-independent.
\end{prop}

\begin{figure}[t]
\begin{tabular}{ccc}
\setlength{\unitlength}{1.2cm}
    \begin{picture}(2,1.7)
      \put(.5,1.015){\framebox(1,.49){0 \boxsp 0}}
      \put(0, .5){\framebox(1,.5){0 \boxsp 1}}
      \put(1.01, .5){\framebox(1,.5){0 \boxsp 0}}
      \put(.5,0){\framebox(1,.485){* \boxsp 0}}
    \end{picture}
& 
\setlength{\unitlength}{1.2cm}
    \begin{picture}(2,1.7)
      \put(.5,1.015){\framebox(1,.49){0 \boxsp 1}}
      \put(0, .5){\framebox(1,.5){1 \boxsp *}}
      \put(1.01, .5){\framebox(1,.5){0 \boxsp 0}}
      \put(.5,0){\framebox(1,.485){1 \boxsp 0}}
    \end{picture}
&
\setlength{\unitlength}{1.2cm}
    \begin{picture}(2,1.7)
      \put(.5,1.015){\framebox(1,.49){0 \boxsp 1}}
      \put(0, .5){\framebox(1,.5){0 \boxsp 0}}
      \put(1.01, .5){\framebox(1,.5){1 \boxsp *}}
      \put(.5,0){\framebox(1,.485){1 \boxsp 0}}
    \end{picture}
\end{tabular}
\caption{The three pairs of the $\pi$-independent Wolfram
  rules that are not $\omega$-independent.\label{fig:minor-forms}}
\end{figure}
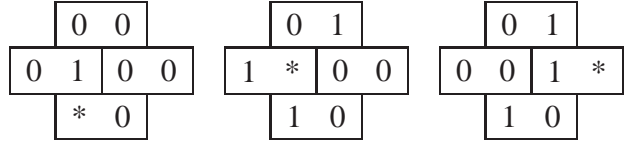

\begin{proof}
Let $k$ be 32 or 40, let $\pi = (\pi_1,\pi_2,\ldots, \pi_n) \in S_n$
be a simple update order, and let $F\colon\F^n\to \F^n$ denote the \sds map
$[\WOLF_n^{(k)},\pi]\colon\F^n\to \F^n$.  The listed rules share the
leftmost form shown in Figure~\ref{fig:minor-forms} and it is easy to
see that $\zer$ is the only fixed state ($\one$ is not fixed and $a_2
= a_6 = 0$ means the rightmost $1$ in any $1$-block converts to $0$
when updated). We also claim $\zer$ is the only periodic state of $F$.
Once this is established, the $\pi$-independence of $\WOLF_n^{(k)}$
follows immediately from Proposition~\ref{prop:fixed}.

The values $a_0 = a_1 = a_4 = 0$ mean non-isolated $0$-blocks persist
indefinitely, they do not shrink or split.  Moreover, $a_2 = a_6 = 0$
means that each non-isolated $0$-block adds at least one $0$ on its
left-hand side with each application of $F$.  In particular, any state
$\y\neq \zer$ with a non-isolated $0$-block eventually becomes the
fixed point $\zer$. Thus no such $\y$ is periodic.

The rest of the argument is by contradiction.  Suppose that $\y$ is a
periodic point of $F$ other than $\zer$ and consider the $i^{\rm th}$
coordinates in $\y$, $F(\y)$ and $F(F(\y))$.  We claim that at least
one of these coordinates is $0$ and at least one of these is $1$.
This is because at least $4$ out of the $5$ local state configurations
that do not involve non-isolated $0$s change the coordinate (and when
$k=32$ all $5$ of them make a change).  The only way that $y_i$ does
not change value in $F(\y)$ is if immediately prior to the application
of $\Wolf_i^{(k)}$, the local state configuration is $011$.  Between
this application of $\Wolf_i^{(k)}$ and the next, the $0$ to the left
is updated.  It either is no longer isolated at this point
(contradicting the periodicity of $\y$) or it now becomes a $1$.  In
the latter case, the application of $\Wolf_i^{(k)}$ during the second
iteration of $F$ changes the $i^{\rm th}$ coordinate from $1$ to $0$.
Note that we used the simplicity of the update order to ensure that
each coordinate is updated only once during each pass through $F$.
Finally, suppose that $i=\pi_1$ and choose $\y$, $F(\y)$ or $F(F(\y))$
so the $(i+1)^{\rm st}$ coordinate is a $0$.  As soon as
$\Wolf_{\pi_1}^{(k)}$ is applied, there is a non-isolated $0$-block,
contradicting the claim that $\y\neq \zer$ is a periodic point.
\end{proof}

Since it was easy to show that every state is periodic under the
bijective Wolfram rule 156 (with tag {\tt -x--}), we did not examine
the evolution of its blocks.  We do so now since its behavior is
relevant to our study of the $4$ remaining rules.

\begin{exmp}[Wolfram rule 156]\label{exmp:156}
Because the symmetric part of rule 156 is {\tt --} no isolated blocks
are ever created or destroyed and thus the number of blocks is
invariant under iteration.  Moreover, the four values $a_1 = a_5 = 0$
and $a_2 = a_3 = 1$ mean that substrings of the form $01$ are fixed
indefinitely, leaving the right end of every $0$-block and the left
end of every $1$-block permanently unchanged.  The other type of
boundary can and does move since $p_3$ = {\tt x}, and it is its
behavior that we want to examine.  Let $\pi \in S_n$ be a simple
update order and let $F\colon\F^n\to \F^n$ denote the \sds map
$[\WOLF_n^{(156)},\pi]\colon\F^n\to \F^n$.  So long as $\y$ is not
$\zer$ or $\one$, there is a $1$-block followed by a $0$-block and a
corresponding substring of the form $01 \cdots 10 \cdots 01$.  (If
$\y$ only contains one $0$-block and one $1$-block, then the first two
digits are the same as the last two digits, but that is irrelevant
here.)  As remarked above, the beginning of the $1$-block and the end
of the $0$-block are fixed, but the other boundary between them can
vary.

Suppose both blocks are non-isolated and consider the central
substring $10$ at positions $i$ and $i+1$.  These are the only
positions in the entire substring that can vary and the first one to
be updated \emph{will} change value.  Assume the $0$ is updated first.
The $1$-block grows, the $0$-block shrinks and the boundary shifts one
step to the right.  As we cycle through the local functions, the
simplicity of $\pi$ guarantees that the $(i+2)^{\rm nd}$ coordinate
is updated before the $(i+1)^{\rm st}$ coordinate is updated a second
time.  Thus the boundary shifts one more step to the right.  This
argument continues to be applicable until the $0$-block shrinks to an
isolated $0$.  At this point, the $0$ is still updated before the $1$
to its left is updated again, but this time the $0$ remains unchanged.
When the $1$ to its left is updated it changes back to a $0$, the
$1$-block shrinks, the $0$-block grows and the boundary shifts to the
left.  The same argument with left and right reversed shows that now
the $0$-block continues to grow until the $1$-block shrinks to an
isolated $1$, at which point the shifting stops and the boundary
starts shifting back in the other direction.
\end{exmp}

\begin{prop}\label{prop:F}
Rules 152 and 184 are $\pi$-independent.
\end{prop}

\begin{proof}
Let $k$ be 152 or 184, let $\pi = (\pi_1,\pi_2,\ldots, \pi_n) \in S_n$
be a simple update order, and let $F\colon\F^n\to \F^n$ denote the \sds map
$[\WOLF_n^{(k)},\pi]\colon\F^n\to \F^n$.  The listed rules share the second
form shown in Figure~\ref{fig:minor-forms} and it is easy to see that
$\zer$ and $\one$ are the only fixed states (since $a_2 = a_6 = 0$
means the rightmost $1$ in any $1$-block converts to $0$ when
updated). We also claim $\zer$ and $\one$ are the only periodic states
of $F$.  Once this is established, the $\pi$-independence of
$\WOLF_n^{(k)}$ follows immediately from Proposition~\ref{prop:fixed}.

Since isolated blocks are never created, the map $\rho$ that sends
$\y$ to the number of blocks it contains is a non-increasing potential
function for $\WOLF_n^{(k)}$.  Moreover, since the only differences
between rule 156 and rules 152 and 184 are that rule 152 removes
isolated $1$-blocks and rule 184 removes both isolated $1$-blocks and
isolated $0$-blocks, the map $F$ agrees with $[\WOLF_n^{(156)},\pi]$
so long as it is not called upon to update an isolated $1$-block (or
an isolated $0$-block when $k$ = 184).  The long-term behavior of rule
156, however, as described in Example~\ref{exmp:156}, shows that under
iteration every $\y$ not equal to $\zer$ or $\one$ eventually updates
such an isolated block, removing it and decreasing $\rho$, thus
showing that such a $\y$ is not periodic.
\end{proof}

Finally, the argument for Wolfram rules 28 and 29 is a combination of
the difficulties found in the proofs of Propositions~\ref{prop:D}
and~\ref{prop:F}.

\begin{prop}\label{prop:G}
Rules 28 and 29 are $\pi$-independent.
\end{prop}

\begin{proof}
Let $k$ be 28 or 29, let $\pi = (\pi_1,\pi_2,\ldots, \pi_n) \in S_n$
be a simple update order, and let $F\colon\F^n\to \F^n$ denote the \sds map
$[\WOLF_n^{(k)},\pi]\colon\F^n\to \F^n$.  The listed rules share the
rightmost form shown in Figure~\ref{fig:minor-forms} and the values
$a_5=0$ and $a_2=1$ mean that isolated blocks are never removed.  Thus
the map $\rho$ that sends $\y$ to the number of blocks it contains is
a non-decreasing potential function for $\WOLF_n^{(k)}$.  The four
values $a_1 = a_5 = 0$ and $a_2 = a_3 = 1$ mean that substrings of the
form $01$ persist indefinitely, as in Wolfram rule 156.  In fact, so
long as $\rho$ is unchanged, the behavior of $F$ under iteration is
indistinguishable from iterations of the map $[\WOLF_n^{(156)},\pi]$.
Consider a substring of the form $01 \cdots 10 \cdots 01$ and suppose
that the length of the $1$-block on the left plus the length of the
$0$-block on the right is at least $4$.  We claim that any $\y$
containing such a substring is not periodic under $F$.  If it were,
the evolution of this substring would oscillate as described in
Example~\ref{exmp:156} and at the point where the $0$-block shrinks to
an isolated $0$, the $1$-block on the left contains the substring
$111$.  Moreover, between the point when that penultimate $0$ becomes
a $1$ and the point when it is to switch back, the substring $111$ is
updated, increasing $\rho$.  When $k$ is 29, a similar increase in
$\rho$ can occur when the $1$-block shrinks to an isolated $1$ and the
$0$-block contains the substring $000$.  In neither case can a state
containing a $1$-block followed by a $0$-block with combined length at
least 4 be periodic under $F$.

Next, note that when $k$ = 29 both of the special states $\zer$ and
$\one$ are not fixed, but that for $k$ = 28 $\one$ is not fixed, while
$\zer$ is fixed.  Let $P$ be the set of states containing both $0$s
and $1$s that do not contain a $1$-block followed by a $0$-block with
combined length at least 4, and, when $k$ = 28, include the special
state $\zer$ as well.  Because we understand the way that such states
$\y\in P$ evolve under Wolfram rule 156 (Example~\ref{exmp:156}), we
know that at no point in the future does a descendent of $\y$ ever
contain a substring of the form $111$ or $000$.  Thus $P$ has been
restricted enough to make the values of $a_7$ and $a_0$ irrelevant,
and $F$ sends $P$ into itself.  Moreover, since $F$ agrees with
$[\WOLF_n^{(156)},\pi]$ on $P$, and this map is injective, $F$ is
injective on $P$, $F$ permutes the states in $P$ and a sufficiently
high power of $F$ is the identity, showing every state in $P$ is
periodic, independent of our choice of $\pi$.
\end{proof}

The rules in this section are the only $\pi$-independent rules of the
$41$ (up to equivalence) that fail to be $\omega$-independent, and
were identified by Collin Bleak and Kevin Ahrendt.

\begin{prop}[\cite{Ahrendt:09}]
  Rules 32, 40, 152, 184, 28, and 29 are not $\omega$-independent.
\end{prop}

\begin{proof}
  For each $k\in\{32,40,152,184,28,29\}$, it suffices to exhibit a
  state $\y\in\F_2^n$ that is \emph{not} periodic under $[\WOLF_n^{(k)},\pi]$
  for any $\pi\in S_n$ but is periodic upder $[\WOLF_n^{(k)},\omega]$
  for some $\omega\in W_n$. In all of these cases, the update order
  $\omega=(1,1,2,2,\dots,n,n)$ will suffice.

  \emph{Case 1} ($k=32$ or $40$): The constant state $\y=\one$ is fixed
  by $[\WOLF_n^{(k)},\omega]$.

  \emph{Case 2} ($k=152$ or $184$): The state $\y=(1,1,0,0,\dots,0)$ is
  fixed by $[\WOLF_n^{(k)},\omega]$.

  \emph{Case 3} ($k=28$ or $29$): If $n$ is even, then the state
  $\y=(1,1,0,0,1,0,1,0,\dots,1,0)$ is fixed by
  $[\WOLF_n^{(k)},\omega]$. If $n$ is odd, then the state
  $(1,1,0,0,1,0,1,0,\dots,1,0,0)$ is fixed by $[\WOLF_n^{(k)},\omega]$.

  In all of these cases, the given state is fixed under $\omega$
  because for each local function $\Wolf_i^{(i)}$ that does not fix $\y$,
  $\Wolf_i^{(i)}\circ\Wolf_i^{(i)}$ fixes $\y$.
\end{proof}

\section{Concluding Remarks}\label{sec:conclusion}

Now that the proof of Theorem~\ref{thm:104} is complete, we pause to
make a few comments about it and the 104 $\pi$-independent Wolfram
rules it identifies.  For each of the $8$ local state configurations,
Wolfram rule $k$ either leaves the central coordinate unchanged or it
``flips'' its value.  The number of local state configurations that
are flipped in this way is strongly correlated with the probability
that a given rule is $\pi$-independent.  See Table~\ref{table:flips}.
The numbers in the third row are the binomial coefficients
$\binom{8}{i}$, since they clearly count the number of Wolfram rules
with exactly $i$ flips.  The key facts illustrated by
Table~\ref{table:flips} are that virtually all of the rules with at
most $2$ flips are $\pi$-independent, the percentage drops off
rapidly between $2$ and $6$ flips, and $\pi$-independence is very
rare among rules with $6$ or more flips.  In fact, all $5$ such rules
are $\pi$-independent because they are bijective.  It would
interesting to know whether this observation can be quantitatively (or
even qualitatively) extended to a rigorous assertion about more
general \sdss.

\begin{table}[t]
\begin{tabular}{|c|ccc|ccc|ccc|}
\hline
Number of flips & 0 & 1 & 2 & 3 & 4 & 5 & 6 & 7 & 8 \\
\hline
Number of $\pi$-independent rules & 1 & 8 & 26 & 34 & 26 & 4 & 4 & 0 & 1 \\
\hline
Number of rules & 1 & 8 & 28 & 56 & 70 & 56 & 28 & 8 & 1 \\
\hline
Percentage & 100\% & 100\% & 93\% & 61\% & 37\% & 7\% & 14\% & 0\% & 100\%\\
\hline
\end{tabular}
\vspace*{3mm}
\caption{The number of flips and the probability of
  $\pi$-independence.\label{table:flips}}
\end{table}

Next, there are two aspects of Theorem~\ref{thm:104} that we found
slightly surprising.  First, we did not initially expect the set of
rules that were $\omega$- and $\pi$-independent for small values of
$n$ to match exactly the set of rules that were $\omega$- and
$\pi$-independent for all values of $n>3$.  The second surprise was
that the during the course of the proof we found that the Wolfram
rules truly are local rules, in the sense that their set of periodic
points tended to have essentially local characterizations.

Finally, although the focus of this article was solely the
classification of the 104 $\pi$-independent Wolfram rules, and not the
dynamics of these rules per se, many interesting dynamical properties
arose in the course of the proof.  We are currently studying the
dynamics and periodic sets for all 256 Wolfram rules in greater
detail, as well as examining how the sets of periodic states under an
$\pi$-independent Wolfram rule get permuted as the update order is
altered.  The latter situation involves an object called the
\emph{dynamics group} of an $\pi$-independent \sdss.  In a follow-up
paper~\cite{Macauley:11a}, we classify the periodic states of all $86$
$\pi$-independent rules and describe their dynamics groups.


\medskip\noindent
{\bfseries Acknowledgments}

The second author gratefully acknowledges the support of the National
Science Foundation. The first and third authors thank the Network
Dynamics and Simulation Science Laboratory (NDSSL) at Virginia Tech
for the support of this research. We thank Kevin Ahrendt and Collin
Bleak for pointing out that $\pi$-independence does not imply
$\omega$-independence. We had originally misstated a result from
another paper that we thought implied that, and corrected it in this
version.

\providecommand{\bysame}{\leavevmode\hbox to3em{\hrulefill}\thinspace}
\providecommand{\MR}{\relax\ifhmode\unskip\space\fi MR }
\providecommand{\MRhref}[2]{%
  \href{http://www.ams.org/mathscinet-getitem?mr=#1}{#2}
}
\providecommand{\href}[2]{#2}


\end{document}